\documentclass[11pt,a4paper]{article}
\usepackage[latin1]{inputenc}
\usepackage[american]{babel}      
\usepackage[T1]{fontenc}
\usepackage{amsmath,amssymb,amsthm}
\usepackage{graphicx}
\usepackage{cite}
\usepackage{xcolor}   
\usepackage{geometry}
\title{{\bf Semi-discrete moduli of smoothness and their applications in one- and two- sided error estimates}}
         
\author{ {\bf Danilo Costarelli} \quad and \quad {\bf Donato Lavella}\\  
Department of Mathematics and Computer Science \\
            University of Perugia\\
        1, Via Vanvitelli, 06123 Perugia, Italy    \\  
{\small {\tt danilo.costarelli@unipg.it} \quad  - \quad {\tt donato.lavella@studenti.unipg.it}} }

\date{}

\newcommand{\N}{\mathbb{N}}
\newcommand{\R}{\mathbb{R}}

\newcommand{\be}{\begin{equation}}
\newcommand{\ee}{\end{equation}}

\newtheorem{definition}{Definition}[section]
\newtheorem{remark}[definition]{Remark}
\newtheorem{theorem}[definition]{Theorem}

\newtheorem{example}[definition]{Example}
\newtheorem{corollary}[definition]{Corollary}

\begin{document}
\maketitle 

\begin{abstract}
In this paper, we introduce a new semi-discrete modulus of smoothness, which generalizes the definition given by Kolomoitsev and Lomako (KL) in 2023 (in the paper published in the J. Approx. Theory), and we establish very general one- and two- sided error estimates under non-restrictive assumptions for pointwise linear operators.
The proposed results have been proved exploiting the regularization and approximation properties of certain Steklov integrals introduced by Sendov and Popov in 1983.
By the definition of semi-discrete moduli of smoothness here proposed, we derive sharper estimates than those that can be achieved by the classical averaged moduli of smoothness ($\tau$-moduli).
Furthermore, a Rathore-type theorem is established, and a new notion of K-functional is also introduced showing its equivalence with the semi-discrete modulus of smoothness and its realization. One-sided estimates of approximation can be established for classical operators on bounded domains, such as the Bernstein polynomials. In the case of approximation operators on the whole real line, one-sided estimates can be achieved, e.g., for the Shannon sampling (cardinal) series, as well as for the so-called generalized sampling operators. 
\end{abstract}

\section{Introduction} \label{sec1}

Moduli of smoothness are classical tools in approximation theory, widely used to obtain sharp error estimates for various operators in different functional spaces. \\
Notable examples, besides the classical integral moduli, include the $\tau$-moduli, introduced by Sendov and Popov in \cite{SendovPopov}, and the Ditzian-Totik moduli of smoothness, studied in \cite{DitzianTotik}, which are particularly effective in the context of weighted approximation. \\
However, the theory of moduli of smoothness is an active and evolving area of research, and new moduli are regularly introduced in the literature (see, for example, \cite{akgun,kopo,Jafarov,kol,run}) to address specific problems of approximation. \\
A recent development in this direction is the semi-discrete modulus of smoothness introduced by Y. Kolomoitsev and T. Lomako in \cite{KL} (and recently extended in \cite{K2026}), defined for $r, s \in \mathbb{N}$, $1 \leq p < \infty$, an appropriate $\gamma \in (0, 1]$, and $f \in L^p(\mathbb{T})$ (1-periodic functions defined on $\mathbb{T} = [0,1)$), as:
\begin{equation} \label{moduloKL}
    \Omega_{r,s}\left(f, \frac{\gamma}{n}, X_n \right)_p = ||f_{\gamma/n,r}-f||_{\ell_p(X_n)} +\omega_s(f, 1/n)_p,
\end{equation}
where $f_{\gamma/n,r}$ is a special Steklov operator (see \cite{KL} p. 7), $X_n = (x_{k,n})_{k=1}^n \subset \mathbb{T}, \;  n \in \mathbb{N},$ is a set of nodes, and $||\cdot||_{\ell_p(X_n)}$ the associated discrete semi-norm (see  \cite{KL} p. 2), i.e.,\\
\begin{equation}   \label{normadiscretaKL}
        ||f||_{\ell_p(X_n)} = \displaystyle{ \left( \frac{1}{n} \sum_{k=1}^n |f(x_{k,n})|^p \right) ^{1/p}}.
    \end{equation}
Unlike conventional measures of smoothness, this modulus contains simultaneously information on the smoothness of $f$ in $L^p$ and on the discrete behaviour of $f$ at the points in $X_n$.

The information on the discrete behaviour is captured by the quantity $||f-f_{\delta, r}||_{\ell_p(X_n)}$.
The continuous analogue of the previous semi-norm has proven particularly effective in establishing sharp $L^p(\mathbb{T)}$ - error estimates in several problems of approximation (see, for example, \cite{Kolomoitsev2012,KolomoitsevTrigub2012,Trigub1980} and \cite{Trigub2013}). 

Under specific assumptions, Kolomoitsev and Lomako obtain two-sided $L^p$ - error estimates for linear operators and determine the exact order of decay of the $L^p$ - error of approximation in terms of $\Omega_{r,s}\left(f, \frac{\gamma}{n}, X_n \right)_p$. \\
However, their results pertain exclusively to linear operators of the trigonometric type and for $1$-periodic functions on $[0,1)$; moreover, the proofs of the considered approximation theorems rely on results concerning the best trigonometric polynomial approximation. 

For instance, they obtain that, for all $1 < p < \infty$ and $\alpha \in (0,2)$, the following equivalence holds:
\begin{equation}
    ||f-L_n(f)||_p = \mathcal{O}(n^{-\alpha}) \iff \Omega_{1,2} \left(f, \frac{1}{n}, X_n \right)_p = \mathcal{O}(n^{-\alpha}),
\end{equation}
where $L_n(f)$ denotes the (trigonometric) Lagrange interpolation polynomials of $f$. 
\vskip0.2cm

In this paper, we propose a modification of $\Omega_{r,s}$ in order to extend its applicability also to not necessarily trigonometric operators, and also to the case of functions belonging to \( \Lambda^p \) (see \cite{BARDARO2006269}) and \( L^p([a,b]) \) as well. 

For completeness, we recall that the spaces $\Lambda^p$ are suitable non-trivial subspaces of $L^p(\mathbb{R})$ (containing for instance the Sobolev spaces $W^r_p$), that have been introduced in \cite{BARDARO2006269} to find a natural setting in which pointwise operators on $\mathbb{R}$ defined by series (such as the Whittaker operators or the generalized sampling series) converge and the corresponding averaged moduli of smoothness are finite (since the latter property is not always true in $L^p(\mathbb{R})$). 
\vskip0.2cm

The key point of the proposed modification of $\Omega_{r,s}$ consists in the replacement in \eqref{moduloKL} of \( f_{\delta,r} \) with \( \tilde{f}_{\delta,r} \), where \( \tilde{f}_{\delta,r} \) denotes the Steklov averages introduced by Sendov and Popov in \cite{SendovPopov}. In this way, we propose a new tool (that we denote by $\tilde{\Omega}_{r,s}$) for the measurement of the degree of approximation of ``pointwise" operators when functions belonging to $L^p$ are considered, which is more general than the one considered in \cite{KL}. 
\vskip0.2cm

Note that the results established in this paper do not seem to be obtainable for $\Omega_{r,s}$ by the same techniques here employed; this is due to the fact that $f_{\delta,r}$ does not have the same {\em regularization power} of $\widetilde{f}$.

Indeed, the Steklov operators considered in \cite{KL} are such that $f_{\delta, r}$ are absolutely continuous for every $r$, when $f \in L^p$, while $\widetilde{f}_{\delta, r} \in W^r_p$, i.e., the regularity of $\widetilde{f}_{\delta, r}$ increases with the parameter $r$.
\vskip0.2cm

In particular, for $\tilde{\Omega}_{r,s}$ we obtained both one- and two- sided error estimates for linear operators, in the case of $f \in \Lambda^p$ and of $f \in L^p([a,b])$.

In the literature (see \cite{SendovPopov}), in these cases (i.e., when functions have not pointwise
meaning), it is considered a modified version of $L^p$, in which the
relation of equivalence almost everywhere is removed; moreover the order of convergence is usually measured with the above-mentioned averaged moduli of smoothness ($\tau$-moduli).

Here, we also prove that the established error estimates are sharper than those that can be achieved by the mentioned $\tau$-moduli. 

Furthermore, a Rathore-type theorem is obtained, and a new notion of K-functional is also introduced, showing its equivalence with $\widetilde{\Omega}_{r,s}$ and its realization. \\
Hence, as a consequence, sharp $L^p - $ error estimates for classical families of ``pointwise" linear operators can be obtained, including both interpolation and approximation ones, such as Bernstein polynomials, Shannon cardinal series and generalized sampling operators. \\
These estimates are sharper than those obtainable through the classical averaged moduli of smoothness ($\tau$-moduli), and established in \cite{SendovPopov, BARDARO2006269}. 
\vskip0.2cm

%


\section{Preliminaries: moduli of smoothness, $\tau$-moduli and Steklov functions}   \label{sec2}

We begin recalling the definition of the classical moduli of smoothness. Below, we will denote by $A$ the sets $[a,b]$ or $\mathbb{R}.$ 
\vskip0.2cm 

The $r$-th difference, $r \in \N$, of a function $f: A \to \mathbb{R}$ is:
\begin{equation*}
    \Delta_{h}^r(f,x) = \displaystyle{\sum_{k = 0}^r \binom{r}{k}(-1)^{r-k}f(x+kh)}, \quad h >0.
\end{equation*}

Given $f: A \to \mathbb{R}$, $\Delta_{h}^r(f,x)$ is well-defined for all $x\in A_{rh}$, where $A_{rh}$ coincides with $A$, except when $A = [a,b]$, for which $A_{rh}$ := $[a, b-rh].$     

\begin{definition}
    The $r$-th modulus of smoothness of $f \in L^p(A)$, if $1 \leq p< \infty,$ or of $f \in C(A)$, (for compact $A$), if $p = \infty$, can be defined as:
\begin{equation*}
     \omega_{r}(f, \delta)_p = \sup _{\substack{0 < h \leq \delta}} ||\Delta_h^r(f,\cdot)||_p(A_{rh}), \qquad \delta > 0.
\end{equation*}    
\end{definition}
When $p = \infty$, the notation $||\cdot||_{\infty}(A_{rh})$ denotes the supremum calculated for $x \in A_{rh}$, and the corresponding moduli of smoothness are simply denoted by $\omega_{r}(f, \delta)$. Whereas, for $1 \leq p < \infty$, the notation $||\cdot||_{p}(A_{rh})$ means that the integral in the definition of p-norm is computed on the set $A_{rh}$, i.e.,  
\begin{equation*}
    ||f||_p(A_{rh}) = \left( \int_{A_{rh}} | f(t) |^p dt\right)^{1/p}.
\end{equation*}

Note that the introduction of the sets $A_{rh}$ is substantially required only in the case $A=[a,b]$ and for non-periodic functions, to get the finite differences $\Delta^r_h(f, \cdot )$ well-defined. 

The main properties of $\omega_r(f, \delta)_p$ are well-known; below we recall the ones that are used throughout the paper: 
 \begin{enumerate}
\item   $\omega_r(f,\delta)_p \leq \delta^r||f^{(r)}||_p \;  \text{ for all} \; f \in W_p^r(A)$, $\delta>0$, where $W_p^r(A)$ denote the usual Sobolev spaces. \label{ew}
\item Setting $\omega_0(f, \delta)_p$ = $||f||_p$, $0 \leq k\leq r,$ we have:
\begin{equation} \label{est}
\omega_r(f, \delta)_p\, \leq\, 2^{r-k}\, \omega_k(f,\delta)_p, \quad 1 \leq p \leq +\infty.
\end{equation}
\end{enumerate}

We now turn our attention to the $\tau$-moduli of smoothness (see \cite{SendovPopov,BARDARO2006269}).
In what follows, we denote by $M(A)$ the space of all measurable and bounded functions on $A$.
\vskip0.2cm

From now on, when functions $f \in L^p(A)$ (or to its sub-spaces) are considered, we always assume that they are defined pointwisely, namely we remove from $L^p(A)$ the relation of equivalence almost everywhere (as usual when we deal with families of pointwise operators, see \cite{SendovPopov, KL}).
\begin{definition}
    Let $f \in M(A)$ be fixed. 
    The local modulus of smoothness of order $k$ of $f$ at a point $x \in A$ is defined as
    $$
        \omega_k(f, x; \delta)\, :=\, \sup \left\{ |\Delta_h^k f(t)| : t, t + k h \in \left[x - \frac{k \delta}{2},\, x + \frac{k \delta}{2} \right] \cap A \right\},
    $$
    where $\delta \in [0, \frac{b - a}{k}]$ for $A = [a,b]$, or $\delta > 0$ when $A = \mathbb{R}$.
\end{definition}
It is clear that, for bounded functions $f$ we have that $\omega_k(f, x; \delta)$ turns out to be well-defined.
\begin{definition}
        Let $f\in M(A)$ and $1 \leq p < +\infty$. 
    The $L^p$-averaged modulus of smoothness of order $k$ (or $\tau$-modulus) of the function $f$ is:
    \begin{align*}
        \tau_k(f; \delta)_p &= ||\omega_k(f,\bullet;\delta)||_p, 
    \end{align*}
    where $\delta \in [0, \frac{b - a}{k}]$ for $A = [a,b]$, or $\delta > 0$ when $A = \mathbb{R}$.
\end{definition}
When functions $f$ belong to $L^p([a,b])$ and are bounded, the moduli $\tau_k(f; \delta)_p$ turn out always finite. While, when $f \in L^p(\mathbb{R)}$, the moduli $\tau_k(f; \delta)_p$ may be infinite.
To avoid this drawback, the notion of the spaces $\Lambda^p$ have been introduced in \cite{BARDARO2006269}. 
\begin{definition} 
A sequence of points $\Sigma := (x_j)_{j \in \mathbb{Z}} \subset \mathbb{R}$ is called an admissible partition of $\mathbb{R}$ if
    \begin{equation*}
        0 < \underline{\Delta} := \inf_{j \in \mathbb{Z}}(x_j-x_{j-1}) \leq \overline{\Delta} := \sup_{j \in \mathbb{Z}}(x_j-x_{j-1}) < +\infty. 
    \end{equation*}

Let $\Sigma = (x_j)_{j \in \mathbb{Z}}$ be an admissible partition of $\mathbb{R}$ and let $\Delta_j := x_{j}-x_{j-1}.$ The discrete $l^p(\Sigma)$-norm of $f: \mathbb{R} \to \mathbb{R}$ is defined, for $1 \leq p < +\infty,$ by
    \begin{equation*}
        ||f||_{l_p(\Sigma)} := \displaystyle{ \left\{ \sum_{j\in \mathbb{Z}}|f(x_j)|^p\Delta_j\right\}}^{1/p}. 
    \end{equation*}

The space $\Lambda^p$ for $1 \leq p < \infty$ is defined by:
    \begin{equation*}
        \Lambda^p = \left\{ f \in M(\mathbb{R}); \; ||f||_{l^p(\Sigma)} < +\infty \mbox{ for each admissible sequence } \Sigma\right\}.
    \end{equation*}
\end{definition}

    Note that despite being referred to as a norm, $||\cdot||_{l_p(\Sigma)}$ is actually a semi-norm on $\Lambda^p$.
\vskip0.2cm

    It is not straightforward to check from the definition whether a function $f$ belongs to $\Lambda^p$. However, it is possible to prove that the following spaces
$$
    \Omega_p := \left\{ f \in M(\mathbb{R});\; |f(t)| \leq g(t),\; t \in \mathbb{R},\; \text{for } g \in L^p(\mathbb{R}), \text{with } g \text{ non-negative}, \right. 
$$
$$
        \left.  \text{even, and non-increasing on } [0, \infty) \right\}, \qquad 1 \leq p < \infty, 
$$
and the classical Sobolev spaces $W^{r}_{p}(\mathbb{R)}, \;  r \in \mathbb{N}, \; 1 \leq p < \infty, $ are linear proper subspaces of $\Lambda^p$, i.e.,
$$
  \Omega_p  \subsetneq \Lambda^p \subsetneq L^p(\mathbb{R)}, \quad \mbox{and} \quad W^r_p(\mathbb{R})
 \subsetneq \Lambda^p \subsetneq L^p(\mathbb{R)}.
$$
Functions $f$ belonging to $\Lambda^p$ are such that $\tau_k(f; \delta)$ is finite, as stated in Proposition 25 of \cite{ButzerStens}, namely:
\begin{equation*}
\tau_k(f, \delta)_p < +\infty, \quad \delta>0, \quad 1 \leq p < \infty, \quad k \in \mathbb{N}.
\end{equation*}
   
    It is well-known that $\tau_k(f, \delta)_p$ satisfy the usual monotonicity, the sub-additivity, as well as some other standard properties (see \cite{SendovPopov}). Among these, we recall the following useful inequality:
\begin{align}     \label{amr7}
\begin{split}
        \tau_k(f; \delta)_p\, &\leq\, c_k\, \delta^k\, ||f^{(k)}||_{p}, \quad f \in W^k_p(A), 
\end{split}
\end{align}
where the positive constant $c_k$ depends only on $k$.
\vskip0.2cm

The $\tau$-modulus can be compared with the classical modulus of smoothness. Indeed:
\begin{theorem}[\cite{ButzerStens}, Proposition 2.7 and \cite{SendovPopov}, Theorem 1.4]  \label{uppsmooth}
        For $f \in M(\mathbb{R}) \cap L^p(\mathbb{R})$, $1 \leq p < \infty$, and $k \in \mathbb{N}$, one has
$$
\omega_k(f; \delta)_p\, \leq\, \tau_k(f;\delta)_p, \quad \delta > 0.
$$
While, for $f \in M([a,b])$ we get:
$$
 \omega_k(f, \delta)_p\, \leq\, \tau_k(f,\delta)_p\, \leq (b-a)^{1/p}\omega_k(f,\delta).
$$
\end{theorem}

We now turn our attention to the following Steklov-type averages, which will be essential in establishing the main results of this paper (see \cite{SendovPopov}, p. 32 and p. 33).
\begin{definition}
\label{steksp}
    Let $f \in L^p(A)$, with $1 \leq p < \infty$. Let $r \in \mathbb{N}$ and assume that $0< \delta \leq (b-a)/r$ if $A = [a,b]$, or $\delta >0$ when $A = \mathbb{R}$. The Steklov averages of $f$ are defined by:
    \begin{equation*}
    \begin{aligned}
        &\tilde{f}_{\delta,r}(x) \\
        & =: (-\delta)^{-r} \displaystyle \int_0^\delta \ldots \int_0^\delta \sum_{m = 1}^r(-1)^{r-m+1} \binom{r}{m}f \left(x+ \frac{m}{r}(t_1+\ldots+t_r)\right)\, dt_1 \ldots dt_r,
    \end{aligned}
    \end{equation*}
in case of $A=\R$, and
$$
    \begin{aligned}
        &\tilde{f}_{\delta,r}(x) \\
        & =: (-\delta)^{-r} \displaystyle \int_0^\delta \ldots \int_0^\delta \sum_{m = 1}^r(-1)^{r-m+1} \binom{r}{r-m}f \left(x+ m\, \theta_x\, \right)\, dt_1 \ldots dt_r,
    \end{aligned}
$$
where:
$$
\theta_x\ :=\ {t_1\, + ... +\, t_r \over r }\, -\, {x-a \over b-a}\, \delta, 
$$
if $A=[a,b]$.
\end{definition}
Note that in the literature it is possible to find several different definitions of the Steklov-averages (see \cite{timan1963approximation,KL}).
\vskip0.2cm

    Next, we present the main properties of $\tilde{f}_{\delta,r}$ (see Theorem 2.5, p. 31, of \cite{SendovPopov}).
\begin{theorem}     \label{425}
        Let $f \in L^p(A)$, with $1 \leq p < \infty$. For each $r \in \N $ and for each $\delta$, $0< \delta \leq (b-a)/r$ if $A = [a,b]$ or $\delta > 0 $ if $A = \mathbb{R}$, the Steklov averages $\tilde{f}_{\delta,r}$ belong to $L^p(A)$ and have the following properties:
        \begin{enumerate}
            \item $||\tilde{f}_{\delta,r}-f||_p \leq c_1\omega_r(f,\delta)_p$, where $c_1 = c_1(r)$  is a positive constant depending only on $r$;
            \item $\tilde{f}_{\delta,r} \in W^{r}_p(A)$ and
            \begin{equation*}
                ||\tilde{f}_{\delta,r} ^{(s)}||_p \leq c_2 \delta^{-s}\omega_s(f, \delta)_p, \quad s= 1, \ldots, r,
            \end{equation*}
            where $c_2 = c_2(r)$ is a positive constant depending only on $r$.
            \item $|f(x) -\tilde{f}_{\delta,r}(x)| \leq \omega_r(f, x; 2\delta).$
        \end{enumerate}
\end{theorem}
\begin{remark} \rm \label{nuova-remark-confronto}
In \cite{KL} the authors considered the following Steklov-type functions:
\begin{equation*}\label{KL-steklov}
 f_{\delta, r}(x)=(-1)^{r+1}\binom{2r}{r}^{-1}\frac{2}{t}\int_{-t/2}^{t/2}\sum_{\nu=0}^{r-1}(-1)^\nu\binom{2r}{\nu}f\left(x+\frac{r-\nu}{r}u\right)du, 
\end{equation*}
where $r\in\N$, $\delta>0$ and $f$ is a locally integrable periodic function. It is clear that in general, when $f \in L^p(A)$ it turns out that $f_{\delta, r}$ is absolutely continuous, for every $r \in \N$, differently from $\tilde{f}_{\delta,r}$ whose regularity increases as the parameter $r$ increases (as showed in 2. of Theorem \ref{425}). However, even if the Steklov-type functions $f_{\delta, r}$ have a lower regularization power than $\tilde{f}_{\delta,r}$, they satisfy:
\be \label{steklov-KL}
||f_{\delta, r} -f||_p \leq C_1\, \omega_{2r}(f, \delta)_p, \quad \delta>0,
\ee
which is sharper than the inequality achieved in 1. of Theorem \ref{425}.
\end{remark}


\section{A new semi-discrete modulus of smoothness: one- and two- sided error estimates}  \label{sec3}

In this section, we propose a modification of the definition of the semi-discrete modulus of smoothness introduced in \cite{KL}, in order to derive both one- and two-sided estimates for linear operators in the case of functions belonging to $\Lambda^p$ and $L^p([a,b])$. 

\begin{definition}     \label{newmoddef}
Let $f \in \Lambda^p$, $1 \leq p < \infty$ and let $\Sigma_W := (j/W)_{j \in \mathbb{Z}}$, for $W > 0$. We define
\begin{equation*}  
 \tilde{\Omega}_{r,s} \left( f, \frac{1}{W}, \Sigma_W \right)_p := ||\tilde{f}_{1/W,r}-f||_{\ell_p(\Sigma_W)}+\omega_s(f, 1/W)_p.
\end{equation*}
\end{definition}

Before introducing the semi-discrete modulus for $f \in L^p([a,b])$, we need to modify \eqref{normadiscretaKL}. 
        \begin{definition}         \label{normalpint}
Let $f \in L^p([a,b])$, $1 \leq p < \infty$ and let $X_n = (x_{k,n})_{k=0}^n \subset  [a,b]$ a set of points. A discrete semi-norm $||f||_{\ell_p(X_n)}$ of $f$ is
            \begin{equation*}
                ||f||_{\ell_p(X_n)} = \displaystyle{ \left(\frac{b-a}{n} \sum_{k=0}^n |f(x_{k,n})|^p\right) ^{1/p}}.
            \end{equation*} 
        \end{definition}
Now, we can give the following.
   \begin{definition}     \label{newmod2def}
             Let $f \in L^p([a,b])$, $1 \leq p < \infty$ and let $X_n = (x_{k,n})_{k = 0}^n \subset [a,b]$ a set of points, with $x_{0,n} < x_{1,n} < \ldots < x_{n,n}$. Suppose further that
             \begin{equation*}
                 \min_k(x_{k+1,n}-x_{k,n}) \geq \gamma/n, \quad \gamma \in (0, 1], 
             \end{equation*}
             where $x_{n+1,n} = (b-a) + x_{1,n}$. We define
         \begin{equation*}          
             \tilde{\Omega}_{r,s} \left( f, \gamma/n, X_n \right)_p = ||\tilde{f}_{\gamma/n,r}-f||_{\ell_p(X_n)}+\omega_s(f, 1/n)_p.
        \end{equation*}
\end{definition}
Note that in both Definition \ref{newmoddef} and Definition \ref{newmod2def} we can easily consider the following modifications:
$$
 \tilde{\Omega}_{r,s} \left( f, \varphi(w), \Sigma_W \right)_p := ||\tilde{f}_{\varphi(W),r}-f||_{\ell_p(\Sigma_W)}+\omega_s(f, \varphi(W))_p,
$$
and
$$        
 \tilde{\Omega}_{r,s} \left( f, \widetilde{\varphi}(n), X_n \right)_p = ||\tilde{f}_{\widetilde{\varphi}(n),r}-f||_{\ell_p(X_n)}+\omega_s(f, \widetilde{\varphi}(n))_p,
$$
where both ${\varphi}(W)$ and $\widetilde{\varphi}(n)$ converge to zero as $W\to +\infty$ and $n \to +\infty$, respectively. In the above cases, the results presented below are still valid.
\vskip0.2cm

As mentioned before, our goal is to establish one- and two-sided estimates for pointwise linear operators in terms of $\tilde{\Omega}_{r,s}$.

In the case $A$ = $\mathbb{R}$, we will consider linear operators of the following form:
\begin{equation*}
    G_W(f)(x) = \sum_{k = -\infty}^{+\infty}f \left( x_{k,W} \right)\varphi_{k,W}(x), \quad W >0,
\end{equation*}
whereas if $A = [a,b]$, $G_n(f)$ will be given by:
\begin{equation*}
    G_n(f)(x) = \sum_{k = 0}^{n}f \left( x_{k,n} \right)\varphi_{k,n}(x), \quad n \in\mathbb{N}.
\end{equation*}

    To establish our main results, we begin by assuming some natural conditions on the linear operators. In the case of \( f \in \Lambda^p \), we will suppose that:
    \begin{equation}     \label{CL1}
        ||G_W(f)||_p \leq K_1||f||_{\ell_p(\Sigma_W)}, \qquad  f \in \Lambda^p, 
    \end{equation}

    \begin{equation}
    \label{CL2}
        K_2||f||_{\ell_p(\Sigma_W)} \leq ||G_W(f)||_p, \qquad  f \in \Lambda^p, 
    \end{equation}

    \begin{equation}
    \label{CL3}
        ||f-G_W(f)||_p \leq K_3W^{-s}||f^{(s)}||_p, \qquad f\in W_p^{s}(\mathbb{R}), 
    \end{equation}

    where $K_1 := K_1(p) > 0$, $K_2 := K_2(p) > 0$ and $K_3 := K_3(s, p) > 0$. \\

    In the case of $f \in L^p([a,b])$ analogous conditions to the previous ones will be required, i.e.,
    \begin{equation}
    \label{CL1I}
        ||G_n(f)||_p \leq K_1||f||_{\ell_p(X_n)}, \qquad  f \in L^p([a,b]), 
    \end{equation}

    \begin{equation}
    \label{CL2I}
        \quad K_2||f||_{\ell_p(X_n)} \leq ||G_n(f)||_p, \qquad  f \in L^p([a,b]) ,
    \end{equation}

    \begin{equation}
    \label{CL3I}
        ||f-G_n(f)||_p \leq K_3n^{-s}||f^{(s)}||_p, \qquad  f\in W_p^{s}([a,b]), 
    \end{equation}

    where $K_1 := K_1(p) > 0$, $K_2 := K_2(p) > 0$ and $K_3 := K_3(s, p) > 0$.
\vskip0.2cm

    At this stage, we can state one of the main results of this section. We begin considering the operators $G_W$.
    \begin{theorem}  \label{CLtheorem}
Let $r \in \mathbb{N}$, $1 \leq p < \infty$ and $f \in \Lambda^p.$ If the linear operators $G_W$ are such that the conditions $\eqref{CL1}$ and $\eqref{CL3}$ are satisfied with $s \leq r,\; s \in \mathbb{N,}$ then 
\begin{equation}   \label{CLT1}
        ||f - G_W(f)||_p \leq K_1 ||\tilde{f}_{1/W,r}-f||_{\ell_p(\Sigma_W)}+C_1\omega_s(f, 1/W)_p.        
\end{equation}
In particular, 
\begin{equation*}
        ||f - G_W(f)||_p \leq \max \left\{ K_1, C_1 \right\} \tilde{\Omega}_{r,s} \left( f, \frac{1}{W}, \Sigma_W \right)_p.
\end{equation*}
    \noindent
    Moreover, if $\eqref{CL2}$ is satisfied, then
\begin{equation}     \label{CLT2}
        K_2 ||\tilde{f}_{1/W,r}-f||_{\ell_p(\Sigma_W)} -C_1\omega_s(f, 1/W)_p \leq ||f-G_W(f)||_p.
\end{equation}
    \noindent
The positive constant $C_1$ does not depend on $f$ and $W$.
\end{theorem}
\begin{proof}
    First, we prove $\eqref{CLT1}$. Applying Theorem \ref{425}, $\eqref{CL3}$ and $\eqref{CL1}$ yields:
\begin{equation*}
\begin{aligned}    
    ||f-G_W(f)||_p &\leq ||f-\tilde{f}_{\delta,r}||_p+||\tilde{f}_{\delta,r}-G_W(\tilde{f}_{\delta,r})||_p + ||G_W(\tilde{f}_{\delta,r})-G_W(f)||_p \\
    &\leq c_1\omega_r(f,\delta)_p+K_3W^{-s}||\tilde{f}_{\delta,r}^{(s)}||_p+K_1||f-\tilde{f}_{\delta,r}||_{\ell_p(\Sigma_W)}.
\end{aligned}    
\end{equation*}    

Using $\eqref{est}$ together with Theorem \ref{425} results in
\begin{equation} \label{appoggioCL1}
\begin{aligned}    
    ||f-G_W(f)||_p &\leq 2^{r-s}c_1\omega_s(f, \delta)_p+K_3W^{-s}c_2\delta^{-s}\omega_s(f,\delta)_p+K_1||f-\tilde{f}_{\delta,r}||_{\ell_p(\Sigma_W)} \\
    &=\left( 2^{r-s}c_1+ \frac{K_3c_2}{(W\delta)^s}\right) \omega_s(f, \delta)_p + K_1||f-\tilde{f}_{\delta,r}||_{\ell_p(\Sigma_W)}. 
\end{aligned}    
\end{equation}    

At this point, it suffices to take $\delta = 1/W$ in $\eqref{appoggioCL1}$; then $\eqref{CLT1}$ follows setting $C_1 = 2^{r-s}c_1+ K_3c_2$. 
\vskip0.2cm

For the second part of the proof, assume that \eqref{CL2} holds. So, recalling that $\tilde{f}_{\delta,r} \in W^r_p(\mathbb{R})$, we have
\begin{equation*}
    \begin{aligned}
        K_2||\tilde{f}_{\delta,r}-f||_{\ell_p(\Sigma_W)} &\leq ||G_W(\tilde{f}_{\delta,r}-f)||_p.
    \end{aligned}
\end{equation*}

Now, using the triangle inequality, Theorem \ref{425}, $\eqref{CL3}$, $\eqref{est}$ and again Theorem \ref{425}, we obtain 
\begin{equation} \label{CLtheoremapp2}
    \begin{aligned}
        K_2||\tilde{f}_{\delta,r}-f||_{\ell_p(\Sigma_W)} &\leq ||f-\tilde{f}_{\delta,r}||_p+||\tilde{f}_{\delta,r}-G_W(\tilde{f}_{\delta,r})||_p + ||f-G_W(f)||_p \\
    &\leq c_1\omega_r(f,\delta)_p+K_3W^{-s}||\tilde{f}_{\delta,r}^{(s)}||_p+||f-G_W(f)||_p \\
        &\leq (2^{r-s}c_1+K_3(\delta W)^{-s}c_2)\omega_s(f,\delta)_p+||f-G_W(f)||_p \\
        &\leq C_1\omega_s(f,\delta)_p+||f-G_W(f)||_p. \\
    \end{aligned}    
\end{equation}
As before, taking $\delta = 1/W$ in $\eqref{CLtheoremapp2}$ we obtain $\eqref{CLT2}$. 
\end{proof}

In case $f$ belongs to $L^p([a,b])$, shortly denoted by $L^p(I)$, an analogous result holds:
    \begin{theorem}   \label{CLteoI}
    Let $r, n \in \mathbb{N}$, $1 \leq p < \infty$ and $f \in L^p(I).$ If the conditions $\eqref{CL1I}$ and $\eqref{CL3I}$ are satisfied with $s \leq r,$ then 
    \begin{equation}     \label{CLT1I}
        ||f - G_n(f)||_p \leq K_1 ||\tilde{f}_{\gamma/n,r}-f||_{\ell_p(X_n)}+C_1\omega_s(f, 1/n)_p.        
    \end{equation}
    In particular, 
    \begin{equation*}
        ||f - G_n(f)||_p \leq \max \left\{ K_1, C_1 \right\} \tilde{\Omega}_{r,s} \left( f, \frac{\gamma}{n}, X_n \right)_p.
    \end{equation*}
    
    Moreover, if $\eqref{CL2I}$ is satisfied, then
    \begin{equation}     \label{CLT2I}
        K_2 ||\tilde{f}_{\gamma/n,r}-f||_{\ell_p(X_n)} -C_1\omega_s(f, 1/n)_p \leq ||f-G_n(f)||_p.
    \end{equation}
  
    The constant $C_1$ does not depend on $f$ and $n$ and is positive.
\end{theorem}
\begin{proof}
    The proof is totally analogous to that of Theorem \ref{CLtheorem}, with the only difference being that we take $\delta = \gamma/n$ at the end of the proof to get the thesis.
\end{proof}
From Theorem $\ref{CLtheorem}$ and Theorem \ref{CLteoI}, we immediately derive the following corollaries. 
\begin{corollary}
     Let $r \in \mathbb{N}$, $1 \leq p < \infty$ and $f \in \Lambda^p$ (or $f \in L^p([a,b])$). If the linear operators $G_W$ (or $G_n(f)$) are such that the conditions $\eqref{CL1}$, $\eqref{CL2}$  and $\eqref{CL3}$ (or $\eqref{CL1I}$, $\eqref{CL2I}$  and $\eqref{CL3I}$) are satisfied with $s \leq r,\; s \in \mathbb{N,}$ then $G_W(f)$ (or $G_n(f)$) converges to $f$ with respect to the $\|\cdot\|_p$ norm if and only if
\begin{equation*}
         ||\tilde{f}_{1/W,r}-f||_{\ell_p(\Sigma_W)} \to 0, \hskip0.2cm \text{ as } \hskip0.2cm W \to \infty \quad (\mbox{or}\ ||\tilde{f}_{\gamma/n,r}-f||_{\ell_p(X_n)} \to 0, \text{ as }\ n \to \infty).
\end{equation*}
\end{corollary}
\begin{corollary}
     Let $r \in \mathbb{N}$, $1 \leq p < \infty$ and $f \in \Lambda^p$ (or $f \in L^p([a,b])$). Assume that the linear operators $G_W$ (or $G_n(f)$) are such that the conditions $\eqref{CL1}$, $\eqref{CL2}$  and $\eqref{CL3}$ (or $\eqref{CL1I}$, $\eqref{CL2I}$  and $\eqref{CL3I}$) are satisfied with $s \leq r,\; s \in \mathbb{N}.$ If
\begin{equation*}
        \omega_s(f, 1/W)_p = o(||\tilde{f}_{1/W,r}-f||_{\ell_p(\Sigma_W)}) \quad ( \mbox{or}\ \omega_s(f, 1/n)_p = o(||\tilde{f}_{\gamma/n,r}-f||_{\ell_p(X_n)})\, ),
\end{equation*}
 then
\begin{equation*}
        ||f-G_W(f)||_p \sim ||\tilde{f}_{1/W,r}-f||_{\ell_p(\Sigma_W)} \quad ( \mbox{or}\ ||f-G_n(f)||_p \sim ||\tilde{f}_{\gamma/n,r}-f||_{\ell_p(X_n)} ).
\end{equation*}  
\end{corollary}
To derive two-sided estimates in terms of $\tilde{\Omega}_{r,s}$, we require an additional result from \cite{RATHORE1994153} (see Theorem 3.1, p.~157), concerning integral functions of finite degree. The class of such functions is denoted by $\mathcal{B}_{\sigma}$ (see \cite{timan1963approximation}, Section 4.3.1, p.~179). We denote by $A_{\sigma}(f)_p$ the error of best approximation of a function $f \in L^p(\mathbb{R})$ by integral functions of finite degree at most $\sigma$, that is,
\[
A_{\sigma}(f)_p = \inf_{g \in \mathcal{B}_\sigma^p} \|f - g\|_p.
\]
Here, $\mathcal{B}^p_{\sigma} \subset \mathcal{B}_{\sigma}$ denotes the set of all complex-valued entire functions of exponential type $\sigma$ whose restriction to $\mathbb{R}$ belongs to $L^p(\mathbb{R})$.
Now, we are able to recall the following result established by Rathore in \cite{RATHORE1994153}, and also known as Rathore theorem.
\begin{theorem} \label{ratlpR}
    Let $f \in L^p(\mathbb{R})$, $1 \leq p < +\infty$, and $r \in \mathbb{N}$. Then:
    \begin{equation}     \label{rathorer1}
        \omega_r ( f, 1/\sigma)_p \leq G A_{\sigma}( f )_p, \quad \sigma >0.
    \end{equation}
    if and only if there exists a positive constant $F$ such that
    \begin{equation}     \label{rathorer2}
        \omega_r ( f, \delta)_p \leq F \omega_{r+1}( f, \delta)_p, \quad \delta > 0.
    \end{equation} 
Similarly, for $f \in L^p(I)$, $1 \leq p \leq  \infty$, and $r \in \mathbb{N}$, we have $E_n(f)_p \sim \omega_r(f, 1/n)_p$ 
if and only if $\omega_r(f, 1/n)_p \sim \omega_{r+1}(f, 1/n)_p$.
\end{theorem}    

Theorem \ref{ratlpR} allows us to prove a Rathore-type theorem for linear operators. 

\begin{theorem}
\label{ratlinopgen}
    Let $r \in \mathbb{N}$, $1 \leq p < \infty$ and $f \in \Lambda^p(\mathbb{R}).$ Assume that the linear operators $G_W$ are such that the conditions $\eqref{CL1}$, $\eqref{CL2}$  and $\eqref{CL3}$ are satisfied with $s \leq r.$ If $G_W(f) \in \mathcal{B}^{p}_W$ and there exists a constant $K >0$ such that
     \begin{equation*}
        \omega_s(f, h)_p \leq  K\omega_{s+1}(f, h)_p, \quad h >0,
     \end{equation*}
     then
     \begin{equation*}
        ||f-G_W(f)||_p \asymp 
        \tilde{\Omega}_{r,s} \left( f, \frac{1}{W}, \Sigma_W \right)_p.
     \end{equation*}
\end{theorem}
\begin{proof}
    The upper bound is a straightforward consequence of $\eqref{CLT1}$.\\ 
    
    As for the lower one, it easily follows from $\eqref{CLT2}$ that
\begin{equation} \label{secondref}
     ||\tilde{f}_{1/W,r}-f||_{\ell_p(\Sigma_W)} + \omega_s(f, 1/W)_p \leq \displaystyle{\frac{1}{K_2} \left( ||f-G_W(f)||_p+(K_2+C_{1}) \omega_s(f, 1/W)_p \right)}.     
\end{equation}

     Then by the first part of Theorem \ref{ratlpR} and using the definition of $A_{\sigma}(f)_p$ it is possible to conclude that:
\begin{equation*}
\begin{aligned}
     ||\tilde{f}_{1/W,r}-f||_{\ell_p(\Sigma_W)} + \omega_s(f, 1/W)_p &\leq \displaystyle{\frac{1}{K_2} \left( ||f-G_W(f)||_p+(K_2+C_{1})GA_W(f)_p\right)}  \\
     &\leq \displaystyle{\frac{1}{K_2} \left( ||f-G_W(f)||_p+(K_2+C_{1})G||f-G_W(f)||_p\right)}  \\
      &\leq \displaystyle{\frac{1+(K_2+C_{1})G}{K_2}||f-G_W(f)||_p}.
\end{aligned}
\end{equation*}
\end{proof}
A Rathore-type theorem for the case of periodic functions belonging to $L^p([0,1])$ has been proved in \cite{KL}; in the case of non-periodic functions defined on the interval $I=[a,b]$ a result similar to Theorem \ref{ratlinopgen} cannot be proved by similar reasoning, since in general a Rathore-type theorem with the error of best approximation by algebraic polynomials and the moduli of smoothness is not valid (some results can be achieved using the Ditzian-Totik moduli of smoothness). 

\vskip0.2cm

Next, we aim to generalize the estimate from above for the semi-discrete modulus of smoothness $\tilde{\Omega}_{r,s}$ in terms of $||f-G_W(f)||_p$ and an appropriate sum, when a certain Bernstein-type inequality is satisfied.
\begin{theorem} \label{newnewtheorem}
    Let $r \in \mathbb{N}$, $1 \leq p < \infty$ and $f \in \Lambda^p.$ Assume that the linear operators $(G_W)_{W > 1}$ are such that the conditions $\eqref{CL1}$, $\eqref{CL2}$ and $\eqref{CL3}$ are satisfied with $s \leq r,\; s \in \mathbb{N}.$ If $G_W(f) \in W^s_p(\mathbb{R})$ and
     \begin{equation}      \label{CL4}
        ||(G_{2^{\nu}}(f) - G_{2^{\nu-1}}(f))^{(s)}||_p \leq K_42^{s\nu}||G_{2^{\nu}}(f)-G_{2^{\nu-1}}(f)||_p, \quad \nu \in \mathbb{N},
     \end{equation}
     with $K_4 = K_4(s,p)$, and $s \leq r$, then
     \begin{equation}      \label{245nuova}
        \begin{aligned}
            &\tilde{\Omega}_{r,s} \left( f, \frac{1}{W}, \Sigma_W \right)_p  \\
            & \leq C \left( ||f-G_W(f)||_p + \displaystyle{\frac{1}{W^s} \sum_{k = 0} ^{[log_2W]}2^{sk}(||f-G_{2^k}(f)||_p + ||f - G_{2^{k-1}}(f)||_p)}\right),
        \end{aligned}
     \end{equation}
     where $C$ is a constant that does not depend on $f$ and $W$.
\end{theorem}
\begin{proof}
    To simplify the notations, denote by $m = [log_2W]$, by $G_{1/2} = 0$ and $G_k = G_k(f), $ $ k \in \mathbb{N}$.

    Firstly, thanks to the subadditivity of the moduli of smoothness and $\eqref{est}$ one has:
    \begin{equation}     \label{2451nuova}
    \begin{aligned}
        \omega_s(f, W^{-1})_p &\leq \omega_s(f-G_{2^m}, W^{-1})_p + \omega_s(G_{2^m} , W^{-1})_p \\
        &\leq 2^s||f-G_{2^m}||_p+\omega_s(G_{2^m}, W^{-1})_p.
    \end{aligned}
    \end{equation}

    Using the technique of the telescoping sums, we can write that
    \begin{equation*}
        G_{2^m} = \sum_{\nu = 0}^{m} (G_{2^\nu} - G_{2^{\nu-1}}). 
    \end{equation*}
    Then, thanks to (1) of Theorem 2.3, $\eqref{ew}$ and by assumption $\eqref{CL4}$, we can estimate $\omega_s(G_{2^m}, W^{-1})_p$ as follows:
    \begin{equation}     \label{2452nuova}
        \begin{aligned}
            \omega_s(G_{2^m}, W^{-1})_p &\leq \sum_{\nu = 0}^{m}\omega_s(G_{2^\nu}-G_{2^{\nu-1}}, W^{-1})_p \\
            &\leq W^{-s}\sum_{\nu = 0}^m||(G_{2^{\nu}}-G_{2^{\nu-1}})^{(s)}||_p \\
            &\leq \displaystyle{\frac{K_4}{W^s}}\sum_{\nu = 0}^m 2^{{s \nu}}||G_{2^{\nu}}-G_{2^{\nu-1}}||_p \\
            &\leq \displaystyle{\frac{K_4}{W^s}}\sum_{\nu = 0}^m 2^{{s \nu}}(||f - G_{2^{\nu}}||_p + ||f-G_{2^{\nu-1}}||_p). \\
       \end{aligned}
    \end{equation}
Putting together $\eqref{2451nuova}$ and $\eqref{2452nuova}$, one obtains that
\begin{equation} \label{246nuova}
\begin{aligned}
    \omega_s(f, W^{-1})_p &\leq 2^s||f - G_{2^m}||_p + \displaystyle{\frac{K_4}{W^s}}\sum_{k = 0}^m 2^{{s k}}(||f - G_{2^{k}}||_p + ||f-G_{2^{k-1}}||_p) \\
    &= 2^s\frac{W^s}{W^s}||f - G_{2^m}||_p + \frac{K_4}{W^s} \sum_{k = 0}^{m} 2^{sk}(||f-G_{2^k}||_p+||f-G_{2^{k-1}}||_p)\\
    &\leq 2^s\frac{(2^{m+1})^s}{W^s}||f - G_{2^m}||_p + \frac{K_4}{W^s} \sum_{k = 0}^{m} 2^{sk}(||f-G_{2^k}||_p+||f-G_{2^{k-1}}||_p)\\
    &= \frac{4^s}{W^s}2^{ms}||f - G_{2^m}||_p + \frac{K_4}{W^s} \sum_{k = 0}^{m} 2^{sk}(||f-G_{2^k}||_p+||f-G_{2^{k-1}}||_p)\\
    &\leq  \frac{K_4+4^s}{W^s} \sum_{k = 0}^{m} 2^{sk}(||f-G_{2^k}||_p+||f-G_{2^{k-1}}||_p),
\end{aligned}
\end{equation}
thus, introducing $\eqref{246nuova}$ into $\eqref{secondref}$, it results:
\begin{equation*}
\begin{aligned}
            &||\tilde{f}_{1/W,r}-f||_{\ell_p(\Sigma_W)} + \omega_s(f, 1/W)_p \\
            & \leq \displaystyle \frac{1}{K_2} ||f-G_W(f)||_p+
            \left( 1+\frac{C_1}{K_2} \right) \frac{K_4+4^s}{W^s} \sum_{k = 0}^{m} 2^{sk}(||f-G_{2^k}||_p+||f-G_{2^{k-1}}||_p)
\end{aligned}
\end{equation*}
and the thesis follows after simple algebraic calculations.
\end{proof}

Also in this case, an analogous result holds if $f \in L^p(I)$. 

\begin{theorem} \label{bernteoI}
    Let $r, n  \; \in \mathbb{N}$, $1 \leq p < \infty$ and $f \in L^p(I).$ Assume that the conditions $\eqref{CL1I}$, $\eqref{CL2I}$  and $\eqref{CL3I}$ are satisfied with $s \leq r.$ If $G_n(f) \in W^s_p(I)$ and
     \begin{equation}      \label{CL4I}
        ||(G_{2^{\nu}}(f) - G_{2^{\nu-1}}(f))^{(s)}||_p \leq K_42^{s\nu}||G_{2^{\nu}}(f)-G_{2^{\nu-1}}(f)||_p, \quad \nu \in \mathbb{N},
     \end{equation}
     with $K_4 = K_4(s,p)$ and $s \leq r$, then
     \begin{equation}      \label{245nuovo}
        \begin{aligned}
            &\tilde{\Omega}_{r,s} \left( f, \frac{\gamma}{n}, X_n \right)_p  \\
            & \leq C \left( ||f-G_n(f)||_p + \displaystyle{\frac{1}{n^s} \sum_{k = 0} ^{[log_2n]}2^{sk}(||f-G_{2^k}(f)||_p + ||f - G_{2^{k-1}}(f)||_p)}\right),
        \end{aligned}
     \end{equation}
     where $C$ is a constant that does not depend on $f$ and $n$.
\end{theorem}

In light of Theorem \ref{newnewtheorem} and Theorem \ref{bernteoI}, the following corollary holds:
\begin{corollary}
    Let $1 \leq p < \infty$, $r\in \mathbb{N}$, $ s \leq r,\; s \in \mathbb{N}$, $\alpha \in (0, s).$ Let $f \in \Lambda^p$ (or $f \in L^p(I)$) and assume that the operators $G_W$, $W > 1$ (or $G_n$, $n \geq 1$), are such that the hypothesis of Theorem \ref{newnewtheorem} (or Theorem \ref{bernteoI}) are satisfied. Then the following assertions are equivalent: \\
    i) $||f-G_W(f)||_p = \mathcal{O}(W^{-\alpha}), \quad W > 1$ \quad (or $||f-G_n(f)||_p = \mathcal{O}(n^{-\alpha})$). \\
    ii) $\displaystyle{\tilde{\Omega}_{r,s} \left( f, \frac{1}{W},\Sigma_W \right)} = \mathcal{O}(W^{-\alpha}), \quad W > 1$ \quad (or $\displaystyle{\tilde{\Omega}_{r,s} \left( f, \frac{\gamma}{n}, X_n \right)_p} = \mathcal{O}(n^{-\alpha})$).   
\end{corollary}


\section{K-functionals: equivalence and realization} \label{sec4}
In this section, we introduce a new semi-discrete $K$-functional and establish its equivalence with $\tilde{\Omega}_{r,s}$, (with realization).

\begin{definition}
\label{funzionaleKnuovanuova}
    Let $f \in \Lambda^p.$ We define
    \begin{equation*}
        \mathcal{K}_s(f,\Sigma_W)_p = \inf_{g\in W_p^s(\mathbb{R})} (||f-g||_{\ell_p(\Sigma_W)}+||f-g||_p+W^{-s}||g^{(s)}||_p).
    \end{equation*}
\end{definition}

Before proving the equivalence between this K-functional and the semi-discrete modulus of smoothness $\tilde{\Omega}_{r,s}$, we recall a useful estimate of $||f||_{\ell_p(\Sigma_W)}$.
\begin{theorem}[\cite{BARDARO2006269}, p. 283, Proposition 14. (iii)]
\label{stimacruciale}
    Let $1 \leq p < \infty$ and $r \in \mathbb{N}$. If $f \in W_p^{r}(\mathbb{R}) \cap C(\mathbb{R})$, then
    \begin{equation*}
        ||f||_{\ell_p(\Sigma_W)} \leq ||f||_{p} + \frac{1}{W}||f'||_{p}.
    \end{equation*}
\end{theorem}

Clearly $\tilde{f}_{\delta,r} \in W^{r}_p(\mathbb{R}) \cap C(\mathbb{R})$ (see \cite{BARDARO2006269}, p. 287, Proposition 21. (i)). 
\vskip0.2cm

Now, the aforementioned equivalence is proved. 
    \begin{theorem}     \label{4312secondref}
        Let $r \in \mathbb{N}$ and  $s \leq r,\; s \in \mathbb{N} $. Let $1 \leq p < \infty$ and $f \in \Lambda^p$. It results
        \begin{equation*}
        \mathcal{K}_s(f,\Sigma_W)_p \asymp \tilde{\Omega}_{r,s} \left( f, \frac{1}{W}, \Sigma_W \right)_p.
        \end{equation*}
    \end{theorem}
    \begin{proof}
        Starting with the upper estimate, we have:
        \begin{equation*}
        \begin{aligned}
            \mathcal{K}_s(f, \Sigma_W)_p &\leq ||f-\tilde{f}_{\delta,r}||_{\ell_p(\Sigma_W)}+ ||f-\tilde{f}_{\delta,r}||_p + W^{-s}||\tilde{f}_{\delta,r}^{(s)}||_p \\
            &\leq ||f-\tilde{f}_{\delta,r}||_{\ell_p(\Sigma_W)}+ c_1\omega_r(f, \delta)_p+(W\delta)^{-s}c_2\omega_s(f, \delta)_p \\
            &\leq ||f-\tilde{f}_{\delta,r}||_{\ell_p(\Sigma_W)}+ (c_12^{r-s}+c_2(W\delta)^{-s})\omega_s(f, \delta)_p.
        \end{aligned}
        \end{equation*}

        Taking $\delta = 1/W$ we obtain the desired estimate. \\

        As for the lower estimate, for every $g \in W^s_p(\mathbb{R)}, $ using Theorem \ref{stimacruciale} and Theorem \ref{425} (iii) one has:
        \begin{equation*}
        \begin{aligned}
            &||\tilde{f}_{1/W,r}-f||_{\ell_p(\Sigma_W)} \\
            &\leq ||\tilde{f}_{1/W,r}-\tilde{g}_{1/W,r}||_{\ell_p(\Sigma_W)}+||\tilde{g}_{1/W,r}-g||_{\ell_p(\Sigma_W)}+||g-f||_{\ell_p(\Sigma_W)} \\
            &= ||\widetilde{(f-g)}_{1/W,r}||_{\ell_p(\Sigma_W)}+ ||\tilde{g}_{1/W,r}-g||_{\ell_p(\Sigma_W)}+||g-f||_{\ell_p(\Sigma_W)} \\
            &\leq ||\widetilde{(f-g)}_{1/W,r}||_p+\frac{1}{W}||\widetilde{(f-g)}'_{1/W,r}||_p + c_3\tau_r \left( g; \frac{r+1}{rW}\right)_p+||f-g||_{\ell_p(\Sigma_W)} \\
            &\leq ||f-g||_p + \frac{1}{W}c_2W\omega_1(f-g, 1/W)_p \\
            &+ c_3\left(2\left(\frac{r+1}{r}+1 \right)\right)^{r+1}\tau_r\left(g, 1/W\right)_p+||f-g||_{\ell_p(\Sigma_W)}. \\
            \end{aligned}
        \end{equation*}

        Now, setting $c = c_3\left(2\left(\frac{r+1}{r}+1 \right)\right)^{r+1}$ and $l = r-s$, we have: 
        \begin{equation*}
        \begin{aligned}
            &||\tilde{f}_{1/W,r}-f||_{\ell_p(\Sigma_W)} \\
            &\leq ||f-g||_p + 2c_2||f-g||_p + c2^l\tau_s\left(g, \frac{r}{r-l}\frac{1}{W}\right)_p+||f-g||_{\ell_p(\Sigma_W)} \\
            &\leq ||f-g||_p + 2c_2||f-g||_p + c2^l\left(2 \left(\frac{r}{r-l}\right)+1\right)^{s+1}\tau_s\left(g,\frac{1}{W}\right)_p+||f-g||_{\ell_p(\Sigma_W)}. \\
        \end{aligned}
        \end{equation*}

        Again, setting $\bar{c} = c2^l\left(2 \left(\frac{r}{r-l}\right)+1\right)^{s+1}$ and using $\eqref{amr7}$ yields
        \begin{equation*}
        \begin{aligned}
            &||\tilde{f}_{1/W,r}-f||_{\ell_p(\Sigma_W)} \\
            &\leq (1+2c_2)||f-g||_p + \bar{c}\tau_s\left(g,\frac{1}{W}\right)_p+||f-g||_{\ell_p(\Sigma_W)}. \\
            &\leq  (1+2c_2)||f-g||_p +c_s\bar{c}W^{-s}||g^{(s)}||_p+||f-g||_{\ell_p(\Sigma_W)}. \\
        \end{aligned}
        \end{equation*}
    Passing to the infimum over all $g \in W_p^s(\mathbb{R})$ we get the thesis.
    \end{proof}
Thanks to the calculations made in the previous theorem, we can now prove that the modulus of smoothness $\tilde{\Omega}_{r,s}$ is sharper than the $\tau$-modulus. 
\begin{theorem}
    \label{stimefini}
    If $s \leq r$, it turns out that
    \begin{equation*}
        \tilde{\Omega}_{r,s}\left( f, \frac{1}{W}, \Sigma_W \right)_p\,  \leq\, C\, \tau_s(f, 1/W)_p.
    \end{equation*}
    \begin{proof}
        Repeating the same steps as in Theorem \ref{4312secondref} and using Theorem \ref{uppsmooth}, it results
        \begin{equation*}
        \begin{aligned} 
            \tilde{\Omega}_{r,s} \left( f, \frac{1}{W}, \Sigma_W \right)_p &= ||\tilde{f}_{1/W,r}-f||_{\ell_p(\Sigma_W)}+\omega_s(f, 1/W)_p \\
            &\leq c_3\tau_r \left( f, \frac{1}{W} +\frac{1}{rW} \right) +  \tau_s(f; 1/W)_p \\
            &\leq \bar{c}\tau_s\left(f,\frac{1}{W}\right)_p + \tau_s(f; 1/W)_p \\
            &= (\bar{c}+1)\tau_s(f; 1/W)_p.
        \end{aligned}
 \end{equation*}
        Setting $C = (\bar{c}+1)$ the thesis follows. 
    \end{proof}
\end{theorem}

In certain applications, it is well-known that moduli of smoothness are equivalent to the so-called realization of the K-functional, namely:
\begin{equation*}
 \omega_s(f, 1/W)_p\ \asymp\ ||f-G_W(f)||_p + W^{-s}||(G_W(f))^{(s)}||_p,
\end{equation*}
where $G_W(f)$ is such that $||f-G_W(f)||_p \underset{\sim}{<} \omega_s(f, 1/W)_p$. This concept generalizes the one introduced for the trigonometric polynomials in \cite{HRIV}.
\vskip0.2cm

In order to obtain the aforesaid result, we will assume that $G_W(f) \in W^s_p(\mathbb{R}),$ and that the following condition holds:
    \begin{equation}     \label{CL5}
        W^{-s}||(G_W(f))^{(s)}||_p\, \leq \,K_5\,\omega_s(G_W(f),W^{-1})_p, \quad f \in \Lambda^p, \quad W > 0,
    \end{equation}
with $K_5 = K_5(s,p)$. 
    
Now we can state the following. 
\begin{theorem}     \label{theorem254bis}
Let $r \in \mathbb{N}$, $s \leq r$, $s \in \mathbb{N}$. Let $1 \leq p < \infty$ and $f \in \Lambda^p$. If the operators $G_W$ satisfy $\eqref{CL1}$, $\eqref{CL2}$, $\eqref{CL3}$, with $s\leq r$, and $\eqref{CL5}$, then
\begin{equation*}   
\tilde{\Omega}_{r,s} \left( f, \frac{1}{W}, \Sigma_W \right)_p  \asymp ||f-G_W(f)||_p + W^{-s}||(G_W(f))^{(s)}||_p.   
\end{equation*}
\end{theorem}
\begin{proof}
From $\eqref{secondref}$, it follows that:
\begin{equation}  \label{app332KLbis}
    \begin{aligned}  
        ||\tilde{f}_{1/W,r}&-f||_{\ell_p(\Sigma_W)} + \omega_s(f, 1/W)_p \\
        &\leq \frac{1}{K_2} \left( ||f-G_W(f)||_p+(K_2+C_{1}) \omega_s(f, 1/W)_p \right) \\  
        &= \frac{1}{K_2} \left( ||f-G_W(f)||_p+(K_2+C_{1}) \omega_s(f+G_W(f)-G_W(f), 1/W)_p \right) \\  
        &\leq \frac{1}{K_2} \left( ||f-G_W(f)||_p+(K_2+C_{1})( \omega_s(f-G_W(f), 1/W)_p + \omega_s(G_W(f), 1/W)_p) \right) \\
        &\leq \frac{1}{K_2} \left( (1+(K_2+C_{1})2^s)||f-G_W(f)||_p+ W^{-s}||(G_W(f))^{(s)}||_p \right). \\
    \end{aligned}  
\end{equation}

This completes the first part of the proof. 

As for the lower bound, using $\eqref{CL5}$, property (1) of Theorem 2.3 and $\eqref{est}$, we have:
    \begin{equation*}
    \begin{aligned}
        W^{-s}||(G_W(f))^{(s)}||_p &\leq K_5\omega_s(G_W(f),W^{-1})_p \\
        &\leq K_5(\omega_s(G_W(f)-f, W^{-1})_p + \omega_s(f, W^{-1})_p) \\
        &\leq K_5(2^s||G_W(f)-f||_p + \omega_s(f, W^{-1})_p).
    \end{aligned}
\end{equation*}

    By Theorem \ref{CLtheorem} and $\eqref{CL5}$, we can deduce that:
    \begin{equation*}
    \begin{aligned}
        ||&G_W(f)-f||_p +W^{-s}||(G_W(f))^{(s)}||_p \\ 
        &\leq (1+2^sK_5)||G_W(f)-f||_p + K_5\omega_s(f, W^{-1})_p \\
        &\leq (1+2^sK_5)(K_1 ||f_{1/W,r}-f||_{\ell_p(\Sigma_W)}+C_1\omega_s(f, 1/W)_p)) + K_5\omega_s(f, W^{-1})_p \\
        &\leq K_1(1+2^sK_5)||f_{1/W,r}-f||_{\ell_p(\Sigma_W)}+(C_1(1+2^sK_5)+ K_5)\omega_s(f, W^{-1})_p, \\
    \end{aligned}
\end{equation*}
    thus proving the estimate from below we were looking for.
\end{proof}

Next, we present a result concerning the estimate of $\tilde{\Omega}_{r,s}$ in terms of an appropriate series. 

\begin{theorem} \label{theorem255nuova}
Let $r \in \mathbb{N}$. Let $1 \leq p < \infty$ and $f \in \Lambda^p$. Assume that the operators $G_W$  satisfy $\eqref{CL1}$, $\eqref{CL2}$, $\eqref{CL3}$, with $s\leq r, \; s \in \mathbb{N}$ and $\eqref{CL5}$. If $G_W(f)$ converges to $f$ in $\| \cdot \|_p$ norm and $G_{W}(f)(\xi) = f(\xi) \text{ for all } \xi \in \Sigma_W$, then
\begin{equation*}   
\tilde{\Omega}_{r,s} \left( f, \frac{1}{W}, \Sigma_W \right)_p  \leq C \sum_{k= 1}^\infty(W2^k)^{-s}||(G_{W2^k}(f))^{(s)}||_p, 
\end{equation*}
with $C$ being a constant independent of $f$ and $W$.
\end{theorem}
\begin{proof}
        Firstly, using $\eqref{CLT2}$ of Theorem \ref{CLtheorem}, one has:
\begin{equation} \label{2551nuova}
    \begin{aligned}
        \|f_{1/W,r}-f\|_{\ell_p(\Sigma_W)} 
        &= \left( \frac{1}{W} \sum_{k=-\infty}^{+\infty} |f_{1/W,r}(k/W)-f(k/W)|^p \right)^{1/p} \\
        &\leq 2^{1/p} \left( \frac{1}{2W} \sum_{k=-\infty}^{+\infty} |f_{1/W,r}(k/2W)-f(k/2W)|^p \right)^{1/p} \\
        &= 2^{1/p}||f_{1/W, r}-f||_{\ell_p(\Sigma_{2W})} \\
        &\leq \frac{2^{1/p}}{K_2}(||f-G_{2W}(f)||_p+C_{1}\omega_s(f, W^{-1})_p).
    \end{aligned}
\end{equation}

Repeating the same steps we did to get to $\eqref{app332KLbis}$, we obtain that:
\begin{equation} \label{2552nuova}
    \omega_s(f, W^{-1})_p \leq 2^s||f-G_{2W}(f)||_p+W^{-s}||(G_{2W}(f))^{(s)}||_p.
\end{equation}

Therefore, putting together $\eqref{2551nuova}$ and $\eqref{2552nuova}$ one has:
\begin{equation} \label{nuovissimissimo}
    \begin{aligned}
        &||f_{1/W, r}-f||_{\ell_p(\Sigma_W)}+\omega_s(f, W^{-1}) \\
        &\leq \frac{2^{1/p}}{K_2}(||f-G_{2W}(f)||_p+C_{1}\omega_s(f, W^{-1})_p) + 2^s||f-G_{2W}(f)||_p+W^{-s}||(G_{2W}(f))^{(s)}||_p \\
        &\leq \left( \frac{2^{1/p}}{K_2}(1+2^sC_{1})+2^s \right)||f-G_{2W}(f)||_p+\left( 1+\frac{2^{1/p}C_{1}}{K_2} \right) W^{-s}||(G_{2W}(f))^{(s)}||_p. \\
    \end{aligned}
\end{equation}

To lighten the notations, we will denote
\begin{equation} \label{Ininuova}
    I_{\nu} = ||G_{2W} - G_{W}(G_{2W})||_p \quad \mbox{ and } \quad G_W = G_{W}(f), \quad W > 0.
\end{equation}

Getting back to assumption $\eqref{CL3}$, it is possible to claim that:
\begin{equation} \label{268nuova}
    I_{W} \leq K_3 W^{-s}||G_{2W}^{(s)}||_p.
\end{equation}

Furthermore, thanks to the fact that that $G_W$ interpolates $f$ on the points $\xi \in \Sigma_W$, one has:
\begin{equation*}
    \begin{aligned}
        (G_{W}f)(x) &= \sum_{k=-\infty}^{\infty}f \left( \frac{k}{W} \right)\varphi_{k,W}(x) \\
        &= \sum_{k=-\infty}^{\infty}(G_{2W}f)\left(\frac{k}{W}\right)\varphi_{k,W}(x) \\
        &= G_{W}(G_{2W}f)(x).
    \end{aligned}
\end{equation*}

Therefore, it is clear that
\begin{equation} \label{ChangeIninuova}
    ||G_{2W}-G_{W}||_p = ||G_{2W}-G_{W}(G_{2W})||_p.
\end{equation}
Using the triangle inequality, one has:

\begin{equation*}
    ||f-G_{W}||_p \leq  ||f-G_{2W}||_p + ||G_{2W} - G_{W}||_p.
\end{equation*}

So, thanks to $\eqref{ChangeIninuova}$, we can conclude that:
\begin{equation*}
    ||f-G_{W}||_p - ||f-G_{2W}||_p  \leq   ||G_{2W} - G_{W}||_p = ||G_{2W}-G_{W}(G_{2W})||_p = I_W.
\end{equation*}

At last, thanks to the fact that the sequence $G_W$ converges to $f$, and using $\eqref{268nuova}$ we get:
\begin{equation*}
\begin{aligned}
    ||f-G_W(f)||_p &= \sum_{k = 0}^{\infty}(||f-G_{W2^k}||_p-||f-G_{W2^{k+1}}||_p) \\
    &\leq \sum_{k=0}^{\infty}I_{W2^k} \leq K_3\sum_{k = 0}^{\infty}(W2^k)^{-s}||G_{W2^{k+1}}^{(s)}||_p \\
    &\leq K_3\sum_{k = 1}^{\infty}(W2^{k-1})^{-s}||G_{W2^{k}}^{(s)}||_p \\
    &\leq K_32^s\sum_{k = 1}^{\infty}(W2^{k})^{-s}||G_{W2^{k}}^{(s)}||_p. \\
\end{aligned}
\end{equation*}

Putting together the above estimate with $\eqref{nuovissimissimo}$ the thesis is then achieved.
\end{proof}
The following corollary holds.
\begin{corollary} \label{corollario-equivalenze-W}
    Let $r \in \mathbb{N}$, $s \leq r$, $s \in \mathbb{N}$, $\alpha \in (0,s)$ and $f \in \Lambda^p$, $1 \leq p < \infty$. If the operators $G_W$ are such that the assumptions made in Theorem \ref{theorem255nuova} are satisfied, the following assertions are equivalent:\\
    i) $||f-G_W(f)||_p = \mathcal{O}(W^{-\alpha}),  \quad W > 1.$ \\
    ii) $\displaystyle{\tilde{\Omega}_{r,s} \left( f, \frac{1}{W}, \Sigma_W \right)_p  = \mathcal{O}(W^{-\alpha}),  \quad W > 1.}$ \\
    iii) $||(G_W(f))^{(s)}||_p = \mathcal{O}(W^{s-\alpha}), \quad W > 1.$
\end{corollary}
\vskip0.2cm

Repeating the proofs of the previous theorems, analogous results to those of this section can be obtained in the case of $f \in L^p(I).$ First, a suitable semi-discrete K-functional is introduced when $f \in L^p(I).$
\begin{definition} \label{funzionaleKI}
    Let $f \in L^p(I).$ We define
\begin{equation*}
        \mathcal{K}_s(f,X_n)_p = \inf_{g\in W_p^s(I)} (||f-g||_{\ell_p(X_n)}+||f-g||_p+n^{-s}||g^{(s)}||_p).
\end{equation*}
\end{definition}

Before proving the equivalence between this K-functional and $\tilde{\Omega}_{r,s}$, we need to establish the following.   
\begin{theorem}
\label{stimacrucialeI}
    Let $1 \leq p < \infty$ and $r \in \mathbb{N}$. If $f \in W^{r}_p(I) \cap C(I)$, and the points in $X_n$ are uniformly distributed in $[a,b]$, with $x_0 = a < x_2 < \ldots < x_n = b$, then
    \begin{equation*}
        ||f||_{\ell_p(X_n)} \leq ||f||_{p} + \left(\frac{b-a}{n}\right)||f'||_{p}.
    \end{equation*}
\end{theorem}
\begin{proof}
    As with Proposition 14, (iii) of \cite{BARDARO2006269}, using the fact that $X_n$ is partition of $[a,b]$ and the mean value theorem yields:
    \begin{equation*}
    \begin{aligned}
    \left( \int_{a}^{b} |f(u)|^p \, du \right)^{1/p}
&=
\left( \sum_{j=1}^{n} \int_{x_{j-1}}^{x_j} |f(u)|^p \, du \right)^{1/p} \\
    &= \left( \frac{b-a}{n}\sum_{j=1}^{n} |f(\xi_j)|^p \right)^{1/p} =: S_1. \\
    \end{aligned}
    \end{equation*}

On the other hand, using Minkowski's and Holder's inequalities with $1/p+1/p' = 1$, it results:
 \begin{equation*}
    \begin{aligned}
    S_2 &:= \left|\left\{ \frac{b-a}{n} \sum_{j = 1}^n|f(x_j)|^p\right\}^{1/p} - \left\{ \frac{b-a}{n}\sum_{j=1}^{n} |f(\xi_j)|^p \right\}^{1/p} \right|  
    \\
    &\leq
\left\{ \frac{b-a}{n}\sum_{j=1}^{n} |f(x_j)-f(\xi_j)|^p \right\}^{1/p} \\
&\leq
\left\{ \frac{b-a}{n}\sum_{j=1}^{n} \left[\int_{x_{j-1}}^{x_j}|f'(u)|du\right]^p \right\}^{1/p} \\
&\leq
\left\{ \frac{b-a}{n}\sum_{j=1}^{n} \int_{x_{j-1}}^{x_j}|f'(u)|^pdu \left( \frac{b-a}{n}\right)^{p/p'}\right\}^{1/p} \leq \frac{b-a}{n}||f'||_p.
    \end{aligned}
    \end{equation*}
\end{proof}
We also need the analogous of $iii)$ of Theorem \ref{425} to the case $[a,b]$.
\begin{theorem}
    Let $f \in L^p ([a,b])$, $1 \leq  p < \infty$. For each integer $r \in \N$ and for each $\delta$, $0< \delta \leq (b-a)/nr$, it results: \\
 \begin{equation*}
        ||\tilde{f}_{\delta,r}-f||_{\ell_p(X_n)} \leq c_3 \tau_r \left(f; \delta+\frac{(b-a)}{rn}\right)_p.
    \end{equation*}
\end{theorem}
\begin{proof}
First, thanks to Theorem 2.5', p.34 of \cite{SendovPopov}, we can write:
\begin{equation*}
\begin{aligned}
    ||f-\tilde{f}_{\delta,r}||_{\ell_p(X_n)} &= \left\{ \frac{b-a}{n}\sum_{j = 1}^n |f(x_j)-\tilde{f}_{\delta,r}(x_j)|^p \right\}^{1/p} \\
    &\leq \left\{ \frac{b-a}{n}\sum_{j = 1}^n [\omega_r(f, x_j; 2\delta)]^p \right\}^{1/p} = ||\omega_r(f, \cdot; 2\delta)||_{\ell_p(X_n)}.
\end{aligned}
\end{equation*}

Following the approach of Proposition 22, p. 289, of \cite{BARDARO2006269}, we have:
\begin{equation*}
\begin{aligned}
    ||\omega_r&(f, \cdot; 2\delta)||_{\ell_p(X_n)} \\
    &\leq  \left\{\sum_{j = 1}^{n}\int_{x_{j-1}}^{x_{j}}[\omega_r(f, x_j; 2\delta)]^pdt \right\}^{1/p} \\
     &\leq  \left\{\sum_{j = 1}^{n}\int_{x_{j-1}}^{x_{j}}\left[\omega_r\left(f, t; 2 \left( \delta+\frac{b-a}{nr} \right)\right) \right]^pdt \right\}^{1/p} \\
     &\leq \left\{\int_{a}^{b}\left[\omega_r\left(f, t; 2 \left( \delta+\frac{b-a}{nr} \right)\right) \right]^pdt \right\}^{1/p} \\
     &= \tau_r \left(f, 2 \left( \delta+\frac{b-a}{nr} \right) \right)_p.
\end{aligned}
\end{equation*}
Finally, the thesis follows by the properties of $\tau_r$. 
\end{proof}
Now, we can immediately state the following.
\begin{theorem} \label{th411-new}
        Let $s \leq r$ and $f \in L^p(I)$, $1 \leq p < +\infty$. It turns out that
$$
        \mathcal{K}_s(f,X_n)_p \asymp \tilde{\Omega}_{r,s} \left( f, \frac{\gamma}{n}, X_n \right)_p.
$$
and there exists $C>0$:
\begin{equation*}
        \tilde{\Omega}_{r,s}\left( f, \frac{\gamma}{n}, X_n \right)_p \leq C\tau_s(f, 1/n)_p.
\end{equation*}
 \end{theorem}
\begin{remark} \rm
Actually, the notion of K-functionals recalled in Definition \ref{funzionaleKI} is the same introduced by Kolomoitsev and Lomako in \cite{KL}. However, the result stated in the first part of Theorem \ref{th411-new} is different from the one established in Theorem 3.4 of \cite{KL}, and this is due to the property (\ref{steklov-KL}) stated in Remark \ref{nuova-remark-confronto}. Indeed, if we consider the case $r=1$ from Theorem \ref{th411-new} we get that:
$$
 \mathcal{K}_1(f,X_n)_p \asymp \tilde{\Omega}_{1,1} \left( f, \frac{\gamma}{n}, X_n \right)_p
$$
while, by means of the moduli of smoothness of Kolomoitsev and Lomako ${\Omega}_{1,s} \left( f, \frac{\gamma}{n}, X_n \right)_p$ the equivalence is established with $\mathcal{K}_s(f,X_n)_p$ for $s=1,2$. Similar considerations can be done for higher valued of $r$. 

Even if the moduli of smoothness of this paper and the ones of \cite{KL} are partially comparable, the achieved estimates for ${\Omega}_{r,s}$ are valid only for trigonometric-type operators (thus for periodic functions) while the estimates achieved with $\widetilde{\Omega}_{r,s}$ are valid for more general operators, not necessarily of the trigonometric type.

  Concerning the case of functions $f$ defined on $\R$ considered in Theorem \ref{4312secondref}, this result was not considered in \cite{KL}, thus (from what we know) it is completely new.

\end{remark}

Also in the case $L^p(I)$, the concept of realization of the K-functional can be discussed. This is when: 
\begin{equation*}
\omega_s(f, 1/n)\,  \asymp\,  ||f-G_n(f)||_p + n^{-s}||(G_n(f))^{(s)}||_p,
\end{equation*}
where $G_n(f)$ is such that $||f-G_n(f)||_p \underset{\sim}{<} \omega_s(f, 1/n)_p.$ \\

In order to establish the realization of the K-functionals, also in this case we assume that $G_n(f) \in W^s_p(I)$ and the condition (similar to $\eqref{CL5}$):
\begin{equation}    \label{CL5I}
        n^{-s}||(G_n(f))^{(s)}||_p \leq K_5\omega_s(G_n(f),n^{-1})_p, \quad f \in L^p(I), \quad n \in \mathbb{N},
    \end{equation}
with $K_5 = K_5(s,p)$, is satisfied. Thus:
\begin{theorem}     \label{teo254newI}
Let $r, n \in \mathbb{N}$. Let $1 \leq p < \infty$ and $f \in L^p(I)$. If $G_n$ and $X_n$ satisfy  $\eqref{CL1I}$, $\eqref{CL2I}$, $\eqref{CL3I}$, with $s\leq r$, and $\eqref{CL5I}$, then
        \begin{equation}         \label{realnew}
        \tilde{\Omega}_{r,s} \left( f, \frac{\gamma}{n}, X_n \right)_p  \asymp ||f-G_n(f)||_p + n^{-s}||(G_n(f))^{(s)}||_p.   
        \end{equation}
\end{theorem}
For the operators $G_n$ also the analogous of Theorem \ref{theorem255nuova} can be proved and then, a corollary analogous to Corollary \ref{corollario-equivalenze-W} could be formulated, but for the sake of brevity, we omit these statements.


\section{Examples and measures of smoothness} \label{sec5}

We now present a broad selection of examples that highlight the fact that estimates by semi-discrete moduli of smoothness are, in general, sharper than those given by $\tau$-moduli. Moreover, also other advantages of $\tilde{\Omega}_{r,s}$ are discussed. 
\vskip0.2cm

Below, when we refer to $\widetilde{f}_{\delta,1}$ we simply write (for brevity) $\widetilde{f}_{\delta}$.
\begin{example} \rm 
The Dirichlet function is
        \[
        f(x) =
        \begin{cases}
        1, & \text{if } x \in \mathbb{Q}, \\
        0, & [a,b] \setminus \mathbb{Q} .
        \end{cases}
        \]

First, note that $||\tilde{f}_{1/n}-f||_{\ell_p(X_n)} = (b-a)^{1/p}$. 
Next, it is easy to check from the definition that:
    \begin{equation*}
        \tau(f,1/n)_p = (b-a)^{1/p}.
    \end{equation*}
    Therefore, 
    \begin{equation*}
        ||\tilde{f}_{1/n}-f||_{\ell_p(X_n)} = \tau(f, 1/n)_p = (b-a)^{1/p}.
    \end{equation*} 
    Furthermore,  $\displaystyle \omega \left(f, \frac{1}{n} \right)_p = 0.$ \\

    So, in this case, it results
    \begin{equation*}
        \tilde{\Omega}_{1,1} \left(f, \frac{1}{n}, X_n \right)_p = ||\tilde{f}_{1/n}-f||_{\ell_p(X_n)} + \omega \left( f, \frac{1}{n} \right)_p = \tau \left( f, \frac{1}{n} \right)_p = (b-a)^{1/p}.
    \end{equation*}
\end{example}

    Next, we study a specific function $f \in L_p([a,b])$ for whom $\tilde{\Omega}_{1,1}\left(f, \frac{1}{2n+1}, X_{2n+1}\right)_p = ||\tilde{f}_{1/(2n+1)}-f||_{\ell_p(X_{2n+1})}+\omega(f, 1/n)_p < \tau(f, 1/n)_p.$
 \begin{example} \rm
 Let
        \[
        f(x) =
        \begin{cases}
        1, & x=\frac{a}{b} \in \mathbb{Q}, \; \gcd(a,b) = 1 \mbox{ and } b \mbox{ is even.} \\
        0, & otherwise  .
        \end{cases}
        \]

    Clearly we have again $\omega(f, 1/n)_p = 0$ since $f \equiv 0$ a.e. in $[a,b]$.
    It turns out that: 
    \begin{equation*}
        ||\tilde{f}_{1/(2n+1)}-f||_{\ell_p(X_{2n+1})} = ||f||_{\ell_p(X_{2n+1})} =  0.
    \end{equation*}

  Finally, from the definition, it easily follows that $\tau(f,1/n)_p = (b-a)^{1/p}.$
 Therefore, it results: 
    \begin{equation*}
    \begin{aligned}
        0 = \Omega_{1,1} \left(f, \frac{1}{2n+1}, X_{2n+1} \right)_p &= ||\tilde{f}_{1/(2n+1)}- f||_{\ell_p(X_{2n+1})}+\omega(f, 1/n)_p \\
        &< \tau(f, 1/n)_p = (b-a)^{1/p}.
    \end{aligned}
    \end{equation*}
\end{example}
The previous examples can be easily extended from $[a,b]$ to the whole $\mathbb{R}$. We now present a last example, dealing with an unbounded function. 
\begin{example} \rm
Let \( 1 \leq p < \infty \), \( 0 < \alpha < 1/p \), \( r \in \mathbb{N} \), and consider the function $f: \mathbb{R} \to \mathbb{R}$, defined as follows: 
\be \label{esempio3}
f(x) =
\begin{cases}
\frac{1}{x^\alpha}, & \text{if } 0 < x < 1, \\
0, & \text{otherwise}.
\end{cases}
\ee
  
First, we want to estimate  the quantity $||f_{1/W}-f||_{\ell_p(\Sigma_W)}$, where $W > 1$. Following the same approach of \cite{KL}, we obtain that
\begin{equation*}
       ||\tilde{f}_{1/W}-f||_{\ell_p(\Sigma_W)} \geq \frac{W^{-(1/p-\alpha)}}{1-\alpha},
\end{equation*}
    and the following upper bound:
    \begin{equation*}
        ||\tilde{f}_{1/W}-f||_{\ell_p(\Sigma_W)} \lesssim   W^{-(1/p- \alpha)}.
    \end{equation*}
   We therefore conclude that:
   \begin{equation*}
        ||\tilde{f}_{1/W}-f||_{\ell_p(\Sigma_W)} \asymp  W^{-(1/p- \alpha)}.
    \end{equation*}
    Then, it is easy to verify that:
    \begin{equation*}
        \omega_r(f,\delta)_p \asymp \delta^{1/p-\alpha}.
    \end{equation*}
    So, in this case
    \begin{equation*}
        \tilde{\Omega}_{1,1} \left(f, \frac{1}{W}, \Sigma_W \right)_p = ||\tilde{f}_{1/W}-f||_{\ell_p(\Sigma_W)} + \omega_s(f,1/W)_p \asymp W^{-(1/p-\alpha)}.
    \end{equation*}
\end{example}

Again, if we consider the function in (\ref{esempio3}) restricted to $[0,1]$ we obtain an example of unbounded function for which $\tilde{\Omega}_{1,1}(f, 1/n, X_n)_p \to 0$, as $n \to +\infty$.

\section{One-sided estimates for approximation operators}   \label{sec6}

Now, we illustrate how the results obtained in the previous sections can be used to study several well-known families of linear operators. 

First, we can immediately observe that for all the examples of trigonometric-type operators presented in \cite{KL} we can still apply all the two-sided estimates established in the Section \ref{sec3} of the present paper. Below, we concentrate our attention only on the new applications, that are mainly related to operators that are not in the trigonometric form.

We begin by considering functions in \(\Lambda^p\) and, in particular, examine the classical Shannon sampling series, i.e.
\begin{equation}         \label{samplingseries}
(S_Wf)(t) := \sum_{k=-\infty}^{\infty} f\left( \frac{k}{W}\right) \operatorname{sinc}(Wt-k), \quad t \in \mathbb{R}, 
\end{equation}
where the well-known $sinc$-function is defined as $\sin(\pi x)/(\pi x)$, for $x \neq 0$, and $1$ for $x=0$.

To apply the results established in Section \ref{sec3}, it is fundamental to prove conditions \eqref{CL1}, \eqref{CL2} and \eqref{CL3}. 

Two of them - that is, \eqref{CL1} and \eqref{CL3} - have been already proved in \cite{BARDARO2006269} (see Corollary 29, p. 294 and Proposition 31, p. 294, respectively). For $f \in \Lambda^p$, $1 < p < \infty$, it turns out that:
$$
            ||S_Wf||_p \leq C||f||_{\ell_p(\Sigma_W)}, \quad
            W > 0,
$$
while for $f \in W^r_p(\mathbb{R}) \cap C(\mathbb{R})$, $r \in \N$, $1<p<+\infty$, there holds
$$
            ||S_Wf-f||_{p} \leq CW^{-r}||f^{(r)}||_p, \quad W > 0, 
$$
where $C$ is a constant not depending on $f$ and $W$. Thus we can immediately state the following.
\begin{theorem}     \label{consequence3}
     Let $r  \in \mathbb{N}$, $1 < p < +\infty$ and $f \in \Lambda^p.$ It results
    \begin{equation*}
        ||f - S_W(f)||_p \leq K_1 ||\tilde{f}_{1/W,r}-f||_{\ell_p(\Sigma_W)}+C_1\, \omega_r(f, 1/W)_p.        
    \end{equation*}
    In particular, 
    \begin{equation*}
        ||f - S_W(f)||_p \leq \max \left\{ K_1, C_1 \right\} \tilde{\Omega}_{r,r} \left( f, \frac{1}{W}, \Sigma_W \right)_p.
    \end{equation*}
    \end{theorem}
   
By Theorem \ref{CLtheorem} the following corollaries hold.
\begin{corollary}
      Let $r  \in \mathbb{N}$, $1 < p < +\infty$ and $f \in \Lambda^p$. If
     \begin{equation*}
         ||\tilde{f}_{1/W,r}-f||_{\ell_p(\Sigma_W)} \to 0, \quad \text{ as } \quad W \to +\infty,
     \end{equation*}
     then
\begin{equation*}
  \lim_{W\to +\infty}  \|S_Wf - f\|_p\ =\ 0.
\end{equation*}    
\end{corollary}
\begin{corollary}
      Let $r  \in \mathbb{N}$, $1 < p < +\infty$ and $f \in \Lambda^p.$ If
     \begin{equation*}
        \omega_r(f, 1/W)_p = \mathcal{O}(||\tilde{f}_{1/W,r}-f||_{\ell_p(\Sigma_W)}), \text{ as }   W \to +\infty,
     \end{equation*}
     then
     \begin{equation*}
        ||f-S_W(f)||_p = \mathcal{O} (||\tilde{f}_{1/W,r}-f||_{\ell_p(\Sigma_W)}), \text{ as }   W \to +\infty.
     \end{equation*}  
\end{corollary}
As an example for the case of $f \in L^p([a, b])$, we will consider the classical Bernstein polynomials $B_n(f)$ on $[0,1]$, given by: 
\[
B_nf(x) = \sum_{k=0}^{n} f\left( \frac{k}{n} \right) \binom{n}{k}x^k(1-x)^n.
\]

Bernstein polynomials are a classical tool in approximation theory and their properties are investigated in many articles (see, for example, \cite{totik1994, chen2017, groot, linsen, ditzian})
Here, we need to prove that conditions \eqref{CL1I} and \eqref{CL2I} are satisfied. From Lemma 4.2., p. 68 of \cite{SendovPopov} we know that, for any bounded function $f$ on $[0,1]$, then
\begin{equation} \label{in-b-1}
            ||B_nf||_p \leq ||f||_{\ell_p(X_n)}.
\end{equation}

As for \eqref{CL3I}, in Lemma 4.3, p. 68 of \cite{SendovPopov} it is proved that, for $f \in W_p^2(I)$, it turns out that:
\begin{equation} \label{in-b-2}
    ||f-B_nf||_p \leq K_3\, n^{-1}\, ||f^{''}||_p.
\end{equation}

Therefore, we can derive the following $L^p$ - error estimates for the Bernstein polynomials $B_nf$. Actually, in order to state the following result the proof of Theorem \ref{CLteoI} must be retraced using the above-mentioned inequalities for the Bernstein polynomials. 
    \begin{theorem}     \label{CLtheoremI}
    Let $1 \leq p < \infty$ and $f \in L^p([0,1]).$ It turns out that 
    \begin{equation*}
        ||f - B_n(f)||_p \leq ||\tilde{f}_{{1 \over \sqrt{n}},2}-f||_{\ell_p(X_n)} + C_1\, \omega_2(f, 1/\sqrt{n})_p,       
    \end{equation*}
with $C_1=K_3+1$, where $K_3$ is the constant of (\ref{in-b-2}).  In particular, 
    \begin{equation*}
        ||f - B_n(f)||_p \leq C_1\, \tilde{\Omega}_{2,2} \left( f, \frac{1}{\sqrt{n}}, X_n \right)_p.
    \end{equation*}
\end{theorem}
\begin{proof}
By standard computations, for $\delta>0$ and using inequalities (\ref{in-b-1}) and (\ref{in-b-2}) together with the properties of the Steklov functions, we get:
$$
\hskip-2cm ||f - B_n(f)||_p \leq ||f - \widetilde{f}_{\delta, 2}||_p + ||\widetilde{f}_{\delta, 2} - B_n(\widetilde{f}_{\delta, 2})||_p + ||B_n(\widetilde{f}_{\delta, 2}) - B_n(f)||_p
$$
$$
\hskip-1.1cm \leq ||f - \widetilde{f}_{\delta, 2}||_p + K_3 n^{-1} \| \widetilde{f}_{\delta, 2}'' \|_p + || f - \widetilde{f}_{\delta, 2}||_{\ell_p(X_n)}
$$
$$
\hskip-0.6cm \leq \omega_2(f, \delta)_p + K_3 n^{-1} \delta^{-2} \omega_2(f,\delta)_p + || f - \widetilde{f}_{\delta, 2}||_{\ell_p(X_n)}.
$$
Thus the proof follows by setting $\delta=1/\sqrt{n}$. 
\end{proof}
Moreover, the following corollaries hold.
\begin{corollary}
     Let $1 \leq p < +\infty$ and $f \in L^p([0,1]).$ If 
\begin{equation*}
  ||\tilde{f}_{1/\sqrt{n},2}-f||_{\ell_p(X_n)} \to 0, \quad \text{ as } \quad n \to \infty.
\end{equation*}
then the sequence $B_n(f)$ converges to $f$ in $L^p([0,1])$.
\end{corollary}
\begin{corollary}
     Let $1 \leq p < +\infty$ and $f \in L^p([0,1])$. If
     \begin{equation*}
        \omega_2(f, 1/\sqrt{n})_p = {\cal O}(||\tilde{f}_{1/\sqrt{n},2}-f||_{\ell_p(X_n)}),
     \end{equation*}
     then
\begin{equation*}
        ||f-B_n(f)||_p\ =\ {\cal O}\left( ||\tilde{f}_{1/\sqrt{n},2}-f||_{\ell_p(X_n)} \right).
\end{equation*}  
\end{corollary}
Now, as a further example for the case of $f \in \Lambda^p$, we will consider the generalized sampling operators $S_W^{\varphi}$, introduced by P. Butzer et al. and then investigated in several works (see, for example, \cite{ries-stens, butzer-stens, butzer-stens93, cantarini-costarelli, corso}). We start with their definition.
    \begin{definition}
    \label{defnucl}
        \begin{enumerate}
            \item We say that $\varphi \in C(\mathbb{R)}$ is a time-limited kernel (for a sampling operator) if 
            \begin{enumerate}
                \item there exist $T_0, T_1 \in \mathbb{R},$ with $T_0 < T_1$, such that $\varphi(t) = 0$ for $t \notin [T_0, T_1]$;
                \item $\varphi$ satisfies the partition of the unity property:
                \begin{equation*} 
                    \sum_{k = - \infty}^{\infty}\varphi(u-k) = 1, \quad u \in \mathbb{R}.
                \end{equation*}
            \end{enumerate}
            \item Let $\varphi$ be a time-limited kernel. Then
            \begin{equation}             \label{samplinggene}
                (S_W^{\varphi}(f))(t) := \sum_{k=-\infty}^{\infty} f\left( \frac{k}{W}\right) \varphi(Wt-k)
            \end{equation}
            are called the generalized sampling operators.
        \end{enumerate}
    \end{definition}

Before presenting applications to the above operators, we recall the following definition. 
\begin{definition}
    Let $\varphi$ be a time-limited kernel. The discrete absolute moment of order 0 of $\varphi$ is
    \begin{equation*}
        m_0(\varphi) : = \sup_{u \in \mathbb{R}}\sum_{k = -\infty}^{\infty}|\varphi(u-k)|.
    \end{equation*}
\end{definition}

It is easy to check that, under the assumption of Definition \ref{defnucl}, we have  $m_0(\varphi) < +\infty$. 
\vskip0.2cm

We recall that inequality $\eqref{CL1}$ was proved for the generalized sampling operators in the following theorem:
\begin{theorem}[\cite{ButzerStens}, p. 42, Proposition 3.2] \label{Bern1}
    Let $1 \leq p < \infty$. We have:
    \begin{equation*}
        ||S_W^{\varphi}f||_p \leq \left\{ {m_0(\varphi)^{1-1/p}||\varphi||^{1/p}_{1}} \right\} ||f||_{\ell_p(\Sigma_W)}, \quad f \in \Lambda^p, \; \; \; W > 0.
    \end{equation*}
\end{theorem}

Moreover, we also need a further assumption on the kernel $\varphi$. We require the so-called vanishing moment condition (also known as Strang-Fix type condition), that is
\begin{equation} \label{vanishingnuc}
    \sum_{k = -\infty}^{\infty}(k-u)^j\varphi(u-k) = 0, \quad j = 1, \ldots, r-1, \; \; u \in \mathbb{R}.
\end{equation}

Now we can state the following result. 

\begin{theorem}[\cite{ButzerStens}, p. 44, Proposition 4.2]
\label{Bern3}
    Let $\varphi$ be a time-limited kernel. Assume that condition $\eqref{vanishingnuc}$ holds for some $r \in \mathbb{N}$. If $1\leq p < \infty$ and $f \in W^r_p(\mathbb{R}),$ then
    \begin{equation*}
        ||S_W^{\varphi}f-f||_p \leq \frac{2m_0(\varphi)T^r}{(r-1)!}W^{-r}||g^{(r)}||_p, \quad W > 0,
    \end{equation*}
    where $T:= \max \left\{ |T_0|, |T_1|\right\}$.
\end{theorem}

Theorem \ref{Bern3} shows that $\eqref{CL3}$ holds. Proving inequality \eqref{CL2} is still an open problem. However, Theorem \ref{Bern1} and Theorem \ref{Bern3} allow us to derive the following estimate.

    \begin{theorem}     \label{consequence3}
     Let $r  \in \mathbb{N}$ for which (\ref{vanishingnuc}) holds, $s \leq r$, $1 \leq p < \infty$ and $f \in \Lambda^p$. It results
    \begin{equation*}
        ||f - S_W^{\varphi}(f)||_p \leq K_1 ||\tilde{f}_{1/W,r}-f||_{\ell_p(\Sigma_W)}+C_1\, \omega_s(f, 1/W)_p.        
    \end{equation*}
    In particular, 
    \begin{equation*}
        ||f - S_W^{\varphi}(f)||_p \leq \max \left\{ K_1, C_1 \right\} \tilde{\Omega}_{r,s} \left( f, \frac{1}{W}, \Sigma_W \right)_p.
    \end{equation*}
    \end{theorem}
   
Similarly to the previous cases, the following corollaries can be stated.

\begin{corollary}
Under the assumptions of Theorem \ref{consequence3}, if
     \begin{equation*}
         ||\tilde{f}_{1/W,r}-f||_{\ell_p(\Sigma_W)} \to 0, \quad \text{ as } \quad W \to +\infty,
     \end{equation*}
     then
      \begin{equation*}
        \lim_{W \to +\infty} \|S^{\varphi}_Wf - f\|_p\ =\ 0.
     \end{equation*}    
\end{corollary}
\begin{corollary}
      Under the assumptions of Theorem \ref{consequence3}, if
     \begin{equation*}
        \omega_s(f, 1/W)_p = \mathcal{O}(||\tilde{f}_{1/W,r}-f||_{\ell_p(\Sigma_W)}), \text{ as }   W \to +\infty,
     \end{equation*}
     then
     \begin{equation*}
        ||f-S_W^{\varphi}(f)||_p = \mathcal{O} (||\tilde{f}_{1/W,r}-f||_{\ell_p(\Sigma_W)}), \text{ as }   W \to +\infty.
     \end{equation*}  
\end{corollary}

\section{Conclusions}  \label{sec7}

    In conclusion, in this work we extended the results obtained by Kolomoitsev and Lomako in \cite{KL} for functions belonging to $L^p(\mathbb{T})$ to the cases of $f \in \Lambda^p$ and of $f \in L^p([a,b])$, and for general linear operators, not necessarily of the trigonometric type.
\vskip0.2cm

    The results here established are strongly based on the better regularization properties of the Steklov averages $\widetilde{f}_{\delta, r}$ than the ones of ${f}_{\delta, r}$ (as discussed in Remark \ref{nuova-remark-confronto}). This fact allows to establish the main results of Section \ref{sec3} and Section \ref{sec4} without using classical results of best approximation for trigonometric-type operators, thus obtaining more general one - and two- sided error estimates. Moreover, also the estimates established the case of approximation operators acting of $f: \R \to \R$ are completely new.
\vskip0.2cm

    In addition, we also proved that the estimates achieved by $\tilde{\Omega}_{r,s}$ are asymptotically sharper than the classical ones, which can be established in terms of the well-known $\tau-$moduli.     
%

\section*{Acknowledgments}

{\small The authors would like to thank the Referee for their very useful suggestions for the improvement of the manuscript. The first author is member of the Gruppo Nazionale per l'Analisi Matematica, la Probabilit\`a e le loro Applicazioni (GNAMPA) of the Istituto Nazionale di Alta Matematica (INdAM), of the network RITA (Research ITalian network on Approximation), and of the UMI (Unione Matematica Italiana) group T.A.A. (Teoria dell'Approssimazione e Applicazioni). 
}

\section*{Funding}

{\small The first author has been partially supported within the (1) "National Innovation Ecosystem grant ECS00000041 - VITALITY", funded by the European Union - NextGenerationEU under the Italian Ministry of University and Research (MUR) and (2) 2025 GNAMPA-INdAM Project "MultiPolExp: Polinomi di tipo esponenziale in assetto multidimensionale e multivoco" (CUP E5324001950001). }

\section*{Conflict of interest/Competing interests}

{\small The authors declare that they have no conflict of interest and competing interest.}

\section*{Availability of data and material and Code availability}

{ \small Not applicable.}

%
\end{document}